\documentclass[a4paper,12pt]{amsart}
\RequirePackage{plautopatch}
\RequirePackage[l2tabu, orthodox]{nag}
\usepackage{bxtexlogo}
\usepackage{amsfonts}
\usepackage{bbm}
\usepackage{amssymb, amsmath,stmaryrd,mathrsfs,mathtools}
\usepackage{tikz,float}
\usepackage{tikz-cd}
\usepackage{amsthm}
\usepackage{enumerate}
\usepackage{enumitem}

\usetikzlibrary{calc}
\usepackage{url}
\usepackage[hidelinks]{hyperref}
\usepackage{graphicx}
\usepackage{xcolor}
\usepackage {diagbox}

\theoremstyle{definition}
\newtheorem{thm}{Theorem}

\newtheorem{lem}[thm]{Lemma}

\newtheorem{dfn}[thm]{Definition}
\newtheorem{rem}[thm]{Remark}

\newtheorem{fac}[thm]{Fact}

\newtheorem*{acknowledgements}{Acknowledgements}

\numberwithin{thm}{section} 
\numberwithin{equation}{section}

\DeclareMathAlphabet{\mymathbb}{U}{BOONDOX-ds}{m}{n}

\newcommand{\oo}{{\omega}^\omega}
\newcommand{\ooo}{[{\omega}]^\omega}
\newcommand{\fin}{[\omega]^{<\omega}}

\newcommand{\on}{\mathpunct{\upharpoonright}}

\newcommand{\bb}{\mathfrak{b}}

\newcommand{\dd}{\mathfrak{d}}
\newcommand{\ee}{\mathfrak{e}}
\newcommand{\pr}{\mathfrak{pr}}

\newcommand{\mm}{\mathfrak{m}}

\newcommand{\covm}{\text{cov}(\mathcal{M})}
\newcommand{\nonm}{\text{non}(\mathcal{M})}
\newcommand{\addm}{\text{add}(\mathcal{M})}

\newcommand{\covn}{\text{cov}(\mathcal{N})}
\newcommand{\nonn}{\text{non}(\mathcal{N})}
\newcommand{\addn}{\text{add}(\mathcal{N})}
\newcommand{\cofn}{\text{cof}(\mathcal{N})}

\newcommand{\R}{\mathbf{R}}

\renewcommand{\lq}{\preceq_T}

\DeclareMathOperator{\add}{add}
\DeclareMathOperator{\non}{non}
\DeclareMathOperator{\cov}{cov}

\newcommand{\none}{\non(\mathcal{E})}

\makeatletter
\renewcommand{\p@enumi}{}
\makeatother

\newcommand{\glcbh}{\mathbf{GLc}(b,h)}

\newcommand{\ilcbh}{\mathbf{ILc}(b,h)}
\newcommand{\llcbh}{\mathbf{LLc}(b,h)}

\newcommand{\sbh}{\mathcal{S}(b,h)}

\newcommand{\bbglc}{\bb^\mathrm{GLc}}
\newcommand{\ddglc}{\dd^\mathrm{GLc}}

\newcommand{\bbllc}{\bb^\mathrm{LLc}}
\newcommand{\ddllc}{\dd^\mathrm{LLc}}
\newcommand{\bbilc}{\bb^\mathrm{ILc}}
\newcommand{\ddilc}{\dd^\mathrm{ILc}}

\newcommand{\lsbh}{\mathcal{LS}(b,h)}

\newcommand{\minllc}{\mathrm{minLLc}}
\newcommand{\supllc}{\mathrm{supLLc}}
\newcommand{\minglc}{\mathrm{minGLc}}
\newcommand{\supglc}{\mathrm{supGLc}}
\newcommand{\minilc}{\mathrm{minILc}}
\newcommand{\supilc}{\mathrm{supILc}}

\newcommand{\gprbh}{\mathbf{GPr}(b,h)}
\newcommand{\iprbh}{\mathbf{IPr}(b,h)}
\newcommand{\lprbh}{\mathbf{LPr}(b,h)}

\newcommand{\gprtrbh}{\mathcal{GPR}(b,h)}
\newcommand{\lprtrbh}{\mathcal{LPR}(b,h)}
\newcommand{\eeg}{\ee^\mathrm{G}}
\newcommand{\prg}{\pr^\mathrm{G}}
\newcommand{\eel}{\ee^\mathrm{L}}
\newcommand{\prl}{\pr^\mathrm{L}}
\newcommand{\eei}{\ee^\mathrm{I}}
\newcommand{\pri}{\pr^\mathrm{I}}

\newcommand{\pdby}{\triangleleft^*}
\newcommand{\ipdby}{\triangleleft^\infty}

\newcommand{\mingpr}{\mathrm{minGPr}}
\newcommand{\supgpr}{\mathrm{supGPr}}
\newcommand{\minlpr}{\mathrm{minLPr}}
\newcommand{\suplpr}{\mathrm{supLPr}}
\newcommand{\minipr}{\mathrm{minIPr}}
\newcommand{\supipr}{\mathrm{supIPr}}

\newcommand{\ip}{\mathbb{IP}}

\newcommand{\seqb}{\operatorname{seq}_{<\omega}(b)}

\newcommand{\lebb}{\mathbb{LE}}

\newcommand{\cal}{\mathcal}

\DeclareMathOperator{\Pred}{Pred}

\title{Notes on slalom prediction}
\author{Takashi Yamazoe}
\address{Graduate School of System Informatics, Kobe University,
	Rokko--dai 1--1, Nada--ku, 657--8501 Kobe, Japan}
\email{212x502x@cloud.kobe-u.jp}
\begin{document}

	\begin{abstract}
		We study a concept of evasion and prediction associated with slaloms, called slalom prediction. This article collects ZFC-provable properties on the slalom prediction.
	\end{abstract}
	\maketitle
	
	\section{Introduction}

	Functions $\omega\to\fin$ are called slaloms and it is known that many cardinal invariants can be characterized using relational systems regarding slaloms (e.g. Theorem \ref{thm_bbglcoh_addn},\ref{thm_ddilc_nonm}). Let us introduce them:

	\begin{dfn}
		
		\begin{itemize}
			\item $\R=\langle X,Y,\sqsubset\rangle$ is a relational system if $X$ and $Y$ are non-empty sets and $\sqsubset \subseteq X\times Y$.
			\item We call an element of $X$ a \text{challenge}, an element of $Y$ a \text{response}, and ``$x\sqsubset y$''  ``$x$ is \text{met by }$y$''.
			\item $F\subseteq X$ is $\R$-unbounded if no response meets all challenges in $F$.
			\item $F\subseteq Y$ is $\R$-dominating if every challenge is met by some response in $F$.
			\item $\R$ is non-trivial if $X$ is $\R$-unbounded and $Y$ is $\R$-dominating. For non-trivial $\R$, define
			\begin{itemize}
				\item $\bb(\R)\coloneq\min\{|F|:F\subseteq X \text{ is }\R\text{-unbounded}\}$, and
				\item $\dd(\R)\coloneq\min\{|F|:F\subseteq Y \text{ is }\R\text{-dominating}\}$.
			\end{itemize}

		\end{itemize}
		
	\end{dfn}

	\begin{dfn}
		Let $b\in(\omega+1)^\omega$ and $h\in\oo$.
		\begin{itemize}
			\item Let $\prod b\coloneq\prod_{n<\omega} b(n)$, $\seqb\coloneq\bigcup_{n<\omega}\prod_{i<n}b(i)$ and $\mathcal{S}(b,h)\coloneq\prod_{n<\omega}[b(n)]^{\leq h(n)}$, the set of all $(b,h)$-slaloms. We often use ``localizer'' instead of ``slalom''. 
			\item For $x\in \prod b$ and $\varphi\in \sbh$, we write:
			\begin{itemize}
				\item $x\in^*\varphi$ if $x(n)\in\varphi(n)$ for all but finitely many $n<\omega$, 
				and
				\item $x\in^\infty\varphi$ if $x(n)\in\varphi(n)$ for infinitely many $n<\omega$.
			\end{itemize}
			\item Denote the relational system $\glcbh\coloneq\langle\prod b,\sbh,\in^* \rangle$ and $\bb^\mathrm{GLc}_{b,h}\coloneq\bb(\glcbh)$ and $\dd^\mathrm{GLc}_{b,h}\coloneq\dd(\glcbh)$. (short for ``global localization'')
			\item Denote the relational system $\ilcbh\coloneq\langle\prod b,\sbh,\in^\infty \rangle$ and $\bb^\mathrm{ILc}_{b,h}\coloneq\bb(\ilcbh)$ and $\dd^\mathrm{ILc}_{b,h}\coloneq\dd(\ilcbh)$. (short for ``infinite localization'')
		\end{itemize}
		To avoid triviality, we always assume $1\leq h(n)<b(n)$ for all $n<\omega$.
	\end{dfn}
	
	There are several types of notation denoting the same cardinal invariants we deal with, such as $\mathfrak{c}^\forall_{b,h}$, $\mathfrak{v}^\exists_{b,h}$, $\dd_h(\in^*)$,  $\bb^\mathrm{Lc}_{b,h}$,  $\dd^\mathrm{aLc}_{b,h}$,  $\mathfrak{sl}_\mathrm{t}(b,h,\mathrm{Fin})$,  $\mathfrak{sl}^\bot_\mathrm{e}(b,h,\mathrm{Fin})$
	(see the list in \cite[Remark 3.3]{CGMRS24}) and our notation is based on Cardona and Mej\'{\i}a's  $\bb^\mathrm{Lc}_{b,h}$, $\dd^\mathrm{aLc}_{b,h}$ in \cite{CM23} (they used the names ``localizations'' and ``anti-localizations'').

	On the other hand, Blass studied the Specker phenomenon in group theory and consequently introduced the combinatorial notion of evasion and prediction in \cite{Bla94} as follows:
	
	\begin{dfn}
		\begin{itemize}
			
			\item A pair $\pi=(D,\{\pi_n:n\in D\})$ is a predictor if $D\in\ooo$ and each $\pi_n$ is a function $\pi_n\colon\omega^n\to\omega$. $\Pred$ denotes the set of all predictors.
			\item $\pi\in \Pred$ predicts $f\in\oo$ if $f(n)=\pi_n(f\on n)$ for all but finitely many $n\in D$. 
			$f$ evades $\pi$ if $\pi$ does not predict $f$.
			\item The prediction number $\pr$ and the evasion number $\ee$ are defined as follows\footnote{While the name ``prediction number'' and the notation ``$\pr$'' are not common, we use them in this article.}:
			\begin{gather*}
				\mathfrak{pr}\coloneq\min\{|\Pi|:\Pi\subseteq \Pred,\forall f\in\oo~\exists\pi\in \Pi~ \pi \text{ predicts } f\},\\
				\mathfrak{e}\coloneq\min\{|F|:F\subseteq\oo,\forall \pi\in \Pred~\exists f\in F~ f\text{ evades }\pi\}.
			\end{gather*}
			
		\end{itemize}
	\end{dfn}
	
	Both localizations and predictions are a kind of a guess, but predictions have the following two differences with localizations:
	\begin{enumerate}
		\item Predictors have their infinite set $D$ on which they guess.
		\item Predictors can see the previous values of functions as a hint of the next value.
	\end{enumerate}

	%

	Now we consider possible patterns of the guess combining the two notions localization and prediction. Blass made a similar approach in \cite{Bla10} (and our names ``global'', ``local'' and ``infinite'' come from his work), but he assumed $b=\omega$ (the constant function of the value $\omega$). Our following framework contains the case when $b\in\oo$ as well.
	
	Firstly, we introduce localizations with their infinite set on which localizers guess:
	
	\begin{dfn}(\cite[Definition 2.4]{BS96}, \cite[Definition 1.2]{Spi98}, \cite[Definition 2.3]{Car24}). 
		Let $b\in(\omega+1)^\omega$ and $h\in\oo$.
		\begin{itemize}
			\item A \text{local $(b,h)$-slalom} is a pair $(D,\varphi)$ where $D\in\ooo$ and $\varphi\in\sbh$. $\lsbh$ denotes the set of all local $(b,h)$-slaloms.
			\item For $x\in \prod b$ and $(D,\varphi)\in\lsbh$, we write $x\in^*(D,\varphi)$  if $x(n)\in\varphi(n)$ for all but finitely many $n\in D$. 
			
			\item Denote the relational system $\llcbh\coloneq\langle\prod b,\lsbh,\in^* \rangle$ and $\bb^\mathrm{LLc}_{b,h}\coloneq\bb(\llcbh)$ and $\dd^\mathrm{LLc}_{b,h}\coloneq\dd(\llcbh)$.
		\end{itemize}
		We often write $x\in^*_D\varphi$ instead of $x\in^*(D,\varphi)$. 
	\end{dfn}
	
	Secondly, we introduce predictions with slaloms and split into three cases (global/local/infinite): 
	\begin{dfn}
		Let $b\in(\omega+1)^\omega$ and $h\in\oo$.
		\begin{itemize}
			\item A global $(b,h)$-predictor is a function $\pi$ with its domain $\seqb$ such that given $n<\omega$ and $\sigma\in\prod_{i<n}b(i)$, $\pi(\sigma)$ is a subset of $b(n)$ of size $\leq h(n)$. $\gprtrbh$ denotes the set of all global $(b,h)$-predictors.
			\item A local $(b,h)$-predictor is a pair $(D,\pi)$ where $D\in\ooo$ and $\pi\in\gprtrbh$. $\mathcal{LPR}(b,h)$ denotes the set of all local $(b,h)$-predictors.
			\item For $(D,\pi)\in\mathcal{LPR}(b,h)$ and $f\in\prod b$, $(D,\pi)$ \text{predicts} $f$ if $f(n)\in \pi(f\on n)$ for all but finitely many $n\in D$ and we write $f\pdby(D,\pi)$ or $f\pdby_D\pi$.
				%
			\item For $\pi\in\mathcal{GPR}(b,h)$ and $f\in\prod b$;
			\begin{itemize}
				\item $\pi$ \text{predicts} $f$ if $f\pdby_\omega\pi$ and we write $f\pdby\pi$.
				\item $\pi$ \text{infinitely predicts} $f$ if $f\pdby_D\pi$ for some $D\in\ooo$ and we write $f\ipdby\pi$. 
			\end{itemize}
			
			\item Define the following relational systems and cardinal invariants:
			\begin{itemize}
				\item $\gprbh\coloneq\langle\prod b,\gprtrbh,\pdby\rangle$, $\eeg_{b,h}\coloneq\bb(\gprbh)$ and $\prg_{b,h}\coloneq\dd(\gprbh)$.
				\item $\lprbh\coloneq\langle\prod b,\lprtrbh,\pdby\rangle$, $\eel_{b,h}\coloneq\bb(\lprbh)$ and $\prl_{b,h}\coloneq\dd(\lprbh)$.
				\item $\iprbh\coloneq\langle\prod b,\gprtrbh,\ipdby\rangle$, $\eei_{b,h}\coloneq\bb(\iprbh)$ and $\pri_{b,h}\coloneq\dd(\iprbh)$.
			\end{itemize}

		\end{itemize}
		
	\end{dfn}
	The local case when $b=\omega$ and $h=1$ is equivalent to the standard prediction. In particular, $\ee=\eel_{\omega,1}$ and $\pr=\prl_{\omega,1}$.
	

	In this article, we focus on ZFC-provable properties on local localizations and slalom predictions, particularly when  $b\in\oo$.
	
	In Section \ref{sec_local}, we study the local localization comparing it with the global and infinite versions following \cite{CM23}. In Section \ref{sec_slalom_predictions}, we see several properties of slalom predictions. In Section \ref{sec_connection}, we investigate the connection between localizations and slalom predictions. In Section \ref{sec_ideals}, we study the relationship between these numbers and some ideals on the reals and conclude the article with Figure \ref{fig_numbers} which illustrates the inequalities we know.
	
	In the rest of this section, we review a basic property of relational systems.

	\begin{dfn}
		
		For relational systems $\R=\langle X,Y,\sqsubset \rangle, \R^{\prime}=\langle X^{\prime},Y^{\prime},\sqsubset^{\prime}~\rangle$,
		$(\Phi_-,\Phi_+):\R\rightarrow\R^\prime$ is a Tukey connection from $\R$ into $\R^{\prime}$ if $\Phi_-:X\rightarrow X^{\prime}$ and $\Phi_+:Y^{\prime}\rightarrow Y$ are functions such that:
		\begin{equation*}
			\forall x\in X~\forall y^{\prime}\in Y^{\prime}~\Phi_-(x)\sqsubset^{\prime} y^{\prime}\Rightarrow x \sqsubset \Phi_{+} (y^{\prime}).
		\end{equation*}

		We write $\R\preceq_T\R^{\prime}$ if there is a Tukey connection from $\R$ into $\R^{\prime}$ and call $\preceq_T$ the Tukey order.
		Note that $\R\preceq_T\R^{\prime}$ implies $\mathfrak{b}(\R^{\prime})\leq\mathfrak{b}(\R)$ and $\mathfrak{d}(\R)\leq\mathfrak{d}(\R^{\prime})$.
		
	\end{dfn}
	

	\section{Local localization} \label{sec_local}
	
	For any function $x\in\prod b$ and $(b,h)$-slalom $\varphi$, $x\in^*\varphi\Rightarrow x\in^*_D\varphi\Rightarrow x\in^\infty\varphi$ holds for any $D\in\ooo$, so: 
	\begin{lem}
		$\ilcbh\lq\llcbh\lq\glcbh$ and hence $\bb^\mathrm{GLc}_{b,h}\leq \bb^\mathrm{LLc}_{b,h}\leq\bb^\mathrm{ILc}_{b,h}$ and $\dd^\mathrm{ILc}_{b,h}\leq \dd^\mathrm{LLc}_{b,h}\leq\dd^\mathrm{GLc}_{b,h}$ for any $b,h$.
	\end{lem} 
	The monotonicity on $b,h$ also holds:
	\begin{lem}
		\label{lem_loc_monotonicity}
		Let $b,b^\prime\in(\omega+1)^\omega$ and $h,h^\prime\in\oo$.
		If $b(n)\leq b^\prime(n)$ and $h^\prime(n)\leq h(n)$ for all but finitely many $n<\omega$, then
		$\mathbf{GLc}(b,h)\lq\mathbf{GLc}(b^\prime,h^\prime)$, $\mathbf{ILc}(b,h)\lq\mathbf{ILc}(b^\prime,h^\prime)$ and $\mathbf{LLc}(b,h)\lq\mathbf{LLc}(b^\prime,h^\prime)$ and hence $\bb^\mathrm{GLc}_{b^\prime,h^\prime}\leq \bb^\mathrm{GLc}_{b,h}$, $\bb^\mathrm{LLc}_{b^\prime,h^\prime}\leq \bb^\mathrm{LLc}_{b,h}$, $\bb^\mathrm{ILc}_{b^\prime,h^\prime}\leq \bb^\mathrm{ILc}_{b,h}$ and $\dd^\mathrm{GLc}_{b,h}\leq \dd^\mathrm{GLc}_{b^\prime,h^\prime}$, $\dd^\mathrm{LLc}_{b,h}\leq \dd^\mathrm{LLc}_{b^\prime,h^\prime}$ and $\dd^\mathrm{ILc}_{b,h}\leq \dd^\mathrm{ILc}_{b^\prime,h^\prime}$.
	\end{lem}
	
	Global and infinite localizations are studied well in \cite{CM23}, so let us look at the local version.  
	The case when $b=\omega$ is studied in \cite{Bla10}. For example: 
	\begin{fac}(\cite{Bla10}, see also \cite[Lemma 2.5]{BS96})
		\label{fac_bllcoh}
		$\bbllc_{\omega,h}=\min\{\ee,\bb\}$ and $\ddllc_{\omega,h}=\max\{\pr,\dd\}$ if $h\in\oo$ goes to infinity\footnote{While Blass proved only the former equation, it is easy to prove the latter from his proof.}.
	\end{fac}
	
	Let us focus on the case $b\in\oo$.
	Following \cite{CM23}, we investigate what kind of properties that hold for global and infinite localizations will also hold for local cases.
	
	First, they proved the following (non-trivial) Tukey relation of global localizations with different parameters:
	\begin{lem}(\cite[Lemma 3.14]{CM23})
		\label{lem_their_3.14}
		Let $b\in\oo$ and $h,h^+\in\oo$. Assume that $\langle I_k\rangle_{k<\omega}$ is an interval partition of $\omega$ satisfying $h(i)\geq h^+(k)$ for all $k<\omega$ and $i\in I_k$ and define $b^+\in\oo$ by $b^+(k)\coloneq\prod_{i\in I_k} b(i)$. Then, $\mathbf{GLc}(b,h)\lq\mathbf{GLc}(b^+,h^+)$.
	\end{lem}
	
	This lemma holds in the local (and infinite) cases as well:
	\begin{lem}
		\label{lem_corresp_3.14}
		Under the assumption of Lemma \ref{lem_their_3.14}, $\mathbf{LLc}(b,h)\lq\mathbf{LLc}(b^+,h^+)$ and $\mathbf{ILc}(b,h)\lq\mathbf{ILc}(b^+,h^+)$.
	\end{lem}

	\begin{proof}
		Their proof of Lemma \ref{lem_their_3.14} gives functions $\Psi_-\colon\prod b\to \prod b^+$ and $\Psi_+\colon\mathcal{S}(b^+,h^+)\to\mathcal{S}(b,h)$ 
		such that for any $x\in\prod b$, $\varphi\in\mathcal{S}(b^+,h^+)$, $k<\omega$ and $i\in I_k$, $\Psi_-(x)(k)\in \varphi(k)$ implies $x(i)\in\Psi_+(\varphi)(i)$, which witness $\mathbf{GLc}(b,h)\lq\mathbf{GLc}(b^+,h^+)$ and $\mathbf{ILc}(b,h)\lq\mathbf{ILc}(b^+,h^+)$ as well. Moreover, for any $D\in\ooo$, $\Psi(x)\in^* (D,\varphi)$ implies $x\in^*(\bigcup_{k\in D}I_k,\Psi_+(\varphi))$, which witnesses $\mathbf{LLc}(b,h)\lq\mathbf{LLc}(b^+,h^+)$.
	\end{proof}
	
	We introduce the limits of localization numbers: 
	\begin{dfn}
		Let $h\in\oo$ go to infinity and $h^\prime\in\oo$.
		\begin{itemize}
			\item $\minglc_h\coloneq\min\{\bbglc_{b,h}:b\in\oo\}$, $\supglc_h\coloneq\sup\{\ddglc_{b,h}:b\in\oo\}$.
			\item $\minllc_h\coloneq\min\{\bbllc_{b,h}:b\in\oo\}$, $\supllc_h\coloneq\sup\{\ddllc_{b,h}:b\in\oo\}$.
			\item $\minilc_{h^\prime}\coloneq\min\{\bbilc_{b,h^\prime}:b\in\oo\}$, $\supilc_{h^\prime}\coloneq\sup\{\ddilc_{b,h^\prime}:b\in\oo\}$.
		\end{itemize}
	\end{dfn}
	
	They proved that the parameter $h$ is irrelevant for the global and infinite numbers:
	
	\begin{thm}(\cite[Theorem 3.19]{CM23})
		\label{thm_GILc_h_indep}
		\begin{itemize}
			\item For $h,h^\prime\in\oo$ going to infinity, $\minglc_h=\minglc_{h^\prime}$ and $\supglc_h=\supglc_{h^\prime}$.
			\item For $h,h^\prime\in\oo$, $\minilc_h=\minilc_{h^\prime}$ and $\supilc_h=\supilc_{h^\prime}$.
		\end{itemize}
		Thus, we omit the subscripts.
	\end{thm}
	
	This holds in the local case as well:
	\begin{thm}
		\label{thm_Lc_h_indep}
		For $h,h^\prime\in\oo$ going to infinity, $\minllc_h=\minllc_{h^\prime}$ and $\supllc_h=\supllc_{h^\prime}$.
		Thus, we omit the subscripts.
	\end{thm}

	\begin{proof}
		We only prove $\minllc_h=\minllc_{h^\prime}$. Let $h^+\in\oo$ such that $h^+(0)=0$ and $h^+\geq^*h,h^\prime$. By Lemma \ref{lem_loc_monotonicity}, $\minllc_{h^\prime},\minllc_{h}\leq\minllc_{h^+}$, so it suffices to show that $\minllc_{h^+}\leq\minllc_{h},\minllc_{h^\prime}$.
		Let $b\in\oo$ be arbitrary.
		Since $h$ and $h^\prime$ go to infinity, we can find some interval partitions $\langle I_k\rangle_{k<\omega}$ and $\langle I_k^\prime\rangle_{k<\omega}$ as in Lemma \ref{lem_corresp_3.14} and corresponding $b^+,(b^\prime)^+$ such that $\mathbf{LLc}(b,h)\lq\mathbf{LLc}(b^+,h^+)$ and $\mathbf{LLc}(b,h^\prime)\lq\mathbf{LLc}((b^\prime)^+,h^+)$. Therefore we have $\minllc_{h^+}\leq\bbglc_{b^+,h^+}\leq\bbglc_{b,h}$ and  $\minllc_{h^+}\leq\bbglc_{(b^\prime)^+,h^+}\leq\bbglc_{b,h^\prime}$.
		Since $b$ was arbitrary, we get $\minllc_{h^+}\leq\minllc_{h},\minllc_{h^\prime}$. 
	\end{proof}
	
	
	They proved that these limits are used to characterize the global localization cardinals for $b=\omega$:
	
	\begin{thm}(\cite[Theorem 3.20]{CM23})
		$\bbglc_{\omega,h}=\min\{\bb,\minglc_h\}$ and $\ddglc_{\omega,h}=\max\{\dd,\supglc_h\}$ for $h\in\oo$ going to infinity. 
	\end{thm}
	
	This theorem holds in the local cases as well:
	\begin{thm}
		\label{thm_b_e_b_minllc}
		$\bbllc_{\omega,h}=\min\{\bb,\minllc_h\}$ and $\ddllc_{\omega,h}=\max\{\dd,\supllc_h\}$ for $h\in\oo$ going to infinity. 
		Equivalently, $\min\{\bb,\ee\}=\min\{\bb,\minllc\}$ and $\max\{\dd,\pr\}=\max\{\dd,\supllc\}$ by Fact \ref{fac_bllcoh}.
	\end{thm}
	
	\begin{proof}
		We only prove $\bbllc_{\omega,h}=\min\{\bb,\minllc_h\}$ and it is clear that $\bbllc_{\omega,h}=\min\{\bb,\ee\}\leq\bb$ and $\bbllc_{\omega,h}\leq \minllc_h$ hold, so we show $\bbllc_{\omega,h}\geq\min\{\bb,\minllc_h\}$. Let $F\subseteq\oo$ of size $<\min\{\bb,\minllc_h\}$. Since $|F|<\bb$, some $b\in\oo$ dominates all $x\in F$ (in the sense of $<^*$) and we may assume that $h\in\prod b$. For $x\in F$, define $x^\prime\in\prod b$ by $x^\prime(n)\coloneq\min\{x(n), b(n)-1\}$ and note that $x(n)=x^\prime(n)$ for all but finitely many $n<\omega$. Since $|F|<\minllc_h\leq\bbllc_{b,h}$, some local $(b,h)$-slalom $(D,\varphi)$ locally localizes all functions in $\{x^\prime:x\in F\}$ and hence all $x\in F$. Therefore, $|F|< \bbllc_{\omega,h}$. 
	\end{proof}
	
	\section{slalom predictions} \label{sec_slalom_predictions}
	For any function $f\in\prod b$ and global $(b,h)$-predictor $\pi$, $f\pdby \pi\Rightarrow f\pdby_D \pi\Rightarrow f\ipdby \pi$ holds for any $D\in\ooo$, so:
	\begin{lem}
		$\iprbh\lq\lprbh\lq\gprbh$ and hence $\eeg_{b,h}\leq \eel_{b,h}\leq\eei_{b,h}$ and $\pri_{b,h}\leq \prl_{b,h}\leq\prg_{b,h}$.
	\end{lem}
	
	The monotonicity on $b,h$ holds as well:
	\begin{lem}
		\label{lem_pre_monotonicity}
		Let $b,b^\prime\in(\omega+1)^\omega$ and $h,h^\prime\in\oo$.
		If $b(n)\leq b^\prime(n)$ and $h^\prime(n)\leq h(n)$ for all but finitely many $n<\omega$, then
		$\mathbf{GPr}(b,h)\lq\mathbf{GPr}(b^\prime,h^\prime)$, $\mathbf{LPr}(b,h)\lq\mathbf{LPr}(b^\prime,h^\prime)$ and $\mathbf{IPr}(b,h)\lq\mathbf{IPr}(b^\prime,h^\prime)$ and hence $\eeg_{b^\prime,h^\prime}\leq \eeg_{b,h}$, $\eel_{b^\prime,h^\prime}\leq \eel_{b,h}$, $\eei_{b^\prime,h^\prime}\leq \eei_{b,h}$ and $\prg_{b,h}\leq \prg_{b^\prime,h^\prime}$, $\prl_{b,h}\leq \prl_{b^\prime,h^\prime}$ and $\pri_{b,h}\leq \pri_{b^\prime,h^\prime}$.
	\end{lem}

	The case when $b=\omega$ is studied in \cite{Bla10}:
	\begin{fac}[\cite{Bla10}]
		\label{fac_Blass_evasion}
		Let $2\leq k<\omega$ and $h\in\oo$ go to infinity.
		\begin{enumerate}
			\item $\aleph_1=\eeg_{\omega,1}\leq\mm_k\leq\eeg_{\omega,k}\leq\addn=\eeg_{\omega,h}$, where $\mm_k$ denotes Martin's number for $\sigma$-$k$-linked forcings. \label{item_Blass_evasion_addn}
			\item $\eel_{\omega,1}=\ee\leq\eel_{\omega,k}\leq\eel_{\omega,h}\leq\covm,\nonm$.
			\item $\eei_{\omega,1}=\eei_{\omega,k}=\eei_{\omega,h}=\covm$. 
		\end{enumerate}
	\end{fac}
	
	It is not hard to see that the dual inequalities hold as follows:
	\begin{lem}
		\label{lem_Blass_evasion}
		Let $2\leq k<\omega$ and $h\in\oo$ go to infinity.
		\begin{enumerate}
			\item $2^{\aleph_0}=\prg_{\omega,1}\geq\prg_{\omega,k}\geq\cofn=\prg_{\omega,h}$. \label{item_Blass_evasion_cofn}
			\item $\prl_{\omega,1}=\pr\geq\prl_{\omega,k}\geq\prl_{\omega,h}\geq\covm,\nonm$.
			\item $\pri_{\omega,1}=\pri_{\omega,k}=\pri_{\omega,h}=\nonm$. 
		\end{enumerate}
	\end{lem}

	The local cases when $b\in\oo$ and $h=1$ are studied in \cite{Bre95}:
	
	\begin{fac}[\cite{Bre95}]
		\label{fac_Brendle_evasion}
		\begin{enumerate}
			\item $\ee\geq\min\{\bb,\ee_{ubd}\}$, where $\ee_{ubd}\coloneqq\min\{\eel_{b,1}:b\in\oo\}$. \label{item_Brendle_evasion_eubd}
			\item All $\eel_{n,1}$ for $2\leq n<\omega$ are the same value $\ee_{fin}$, where the subscript $n$ denotes the constant function of the value $n<\omega$.
			\item $\ee_{fin}\geq\mathfrak{s}$, the splitting number.  
		\end{enumerate}
	\end{fac}
	
	It is not hard to see that the dual inequalities hold as follows:
	\begin{lem}
		\label{lem_Brendle_evasion}
		\begin{enumerate}
			\item $\pr\leq\max\{\dd,\pr_{ubd}\}$, where $\pr_{ubd}\coloneqq\sup\{\prl_{b,1}:b\in\oo\}$. 
			\item All $\prl_{n,1}$ for $2\leq n<\omega$ are the same value $\pr_{fin}$. 
			\item $\pr_{fin}\leq\mathfrak{r}$, the reaping number.  
		\end{enumerate}
	\end{lem}
	
	The same Tukey connections as in Lemma \ref{lem_their_3.14} and \ref{lem_corresp_3.14} also hold by a similar proof:  
	\begin{lem}
		\label{lem_their_3.14_pred}
		Let $b\in\oo$ and $h,h^+\in\oo$. Assume that $\langle I_k\rangle_{k<\omega}$ is an interval partition of $\omega$ satisfying $h(i)\geq h^+(k)$ for all $k<\omega$ and $i\in I_k$ and define $b^+\in\oo$ by $b^+(k)\coloneq\prod_{i\in I_k} b(i)$. Then, $\mathbf{GPr}(b,h)\lq\mathbf{GPr}(b^+,h^+)$, $\mathbf{LPr}(b,h)\lq\mathbf{LPr}(b^+,h^+)$ and $\mathbf{IPr}(b,h)\lq\mathbf{IPr}(b^+,h^+)$.
	\end{lem}
	
	Moreover, in the local and infinite cases, $h$ and $h^\prime$ do not require any assumptions:
	
	\begin{lem}
		\label{lem_their_3.14_pred_1}
		Let $b\in\oo$ and $h,h^+\in\oo$. 
		Then, there exists $b^+\in\oo$ such that $\mathbf{LPr}(b,h)\lq\mathbf{LPr}(b^+,h^+)$ and $\mathbf{IPr}(b,h)\lq\mathbf{IPr}(b^+,h^+)$.
	\end{lem}
	
	\begin{proof}
		We may assume $h=1$. Let $\langle I_k:k<\omega\rangle$ be an interval partition of $\omega$ satisfying $|I_k|>h^+(k)$ for all $k<\omega$ and define $b^+$ by $b^+(k)\coloneq\prod_{i\in I_k}b(i)$.
		Let $\pi$ be any global $(b^+,h^+)$-predictor and fix $k<\omega$ and $\sigma\in \prod_{k^\prime<k}b^+(k^\prime)$.
		Let $S\coloneq\pi(\sigma)\subseteq b^+(k)=\prod_{i\in I_k}b(i)$, which we may assume as a set of functions on $I_k$.
		Since $|S|=|\pi(\sigma)|\leq h^+(k)< |I_k|$, there is $j_k\coloneq j\in I_k$ which is not a branching point of any two $t,t^\prime\in S$, namely, $t \on I_k\cap j=t^\prime\on I_k\cap j$ implies $t(j)=t^\prime(j)$ for any $t,t^\prime\in S$.
		Thus, 
		the function $g_{k,\sigma}:\{t\on I_k\cap j:t\in S\} \to b(j)$, $t\on I_k\cap j\mapsto t(j)$ 
		is well-defined.
		Unfixing $k$ and $\sigma$, let $\pi^\prime$ be a global $(b,1)$-predictor satisfying for all $k<\omega$ and $\tau\in \prod_{i<j_k}b(i)=\prod_{k^\prime<k}b^+(k^\prime)\prod_{i\in I_k\cap j_k}b(i)$, $\pi^\prime(\tau)=\{g_{k,\sigma }(\tau\on I_k\cap j_k)\}$ where $\sigma\coloneqq\langle\tau\on I_{k^\prime}:k^\prime<k\rangle\in \prod_{k^\prime<k}b^+(k^\prime)$.
		By construction, for any $f\in\prod b$ and $k<\omega$,
		\begin{equation*}
			f\on I_k\in\pi(\langle f\on I_{k^\prime}:k^\prime<k\rangle)\text{ implies } f(j_k)\in\pi^\prime (f\on j_k),
		\end{equation*}
		which induces $\mathbf{LPr}(b,1)\lq\mathbf{LLc}(b^+,h^+)$ and $\mathbf{IPr}(b,1)\lq\mathbf{ILc}(b^+,h^+)$. 
	\end{proof}

	We introduce the limits of slalom prediction numbers: 
	\begin{dfn}
		Let $h\in\oo$.
		\begin{itemize}
			\item $\mingpr_h\coloneq\min\{\eeg_{b,h}:b\in\oo\}, \supgpr_h\coloneq\sup\{\prg_{b,h}:b\in\oo\}$. 
			\item $\minlpr_h\coloneq\min\{\eel_{b,h}:b\in\oo\}, \suplpr_h\coloneq\sup\{\prl_{b,h}:b\in\oo\}$.
			\item $\minipr_h\coloneq\min\{\eei_{b,h}:b\in\oo\}, \supipr_h\coloneq\sup\{\pri_{b,h}:b\in\oo\}$. 
		\end{itemize}
		
	\end{dfn}
	%
	%

	By a similar proof to Theorem \ref{thm_b_e_b_minllc}, we obtain the following lemma:

	\begin{lem}
		Let $h\in\oo$.
		\begin{itemize}
			\item $\eeg_{\omega,h}=\min\{\bb,\mingpr_h\}$ and $\prg_{\omega,h}=\max\{\dd,\supgpr_h\}$.
			\item $\eel_{\omega,h}\geq\min\{\bb,\minlpr_h\}$ and $\prl_{\omega,h}\leq\min\{\dd,\suplpr_h\}$. 
			\item $\covm=\eei_{\omega,h}\geq\min\{\bb,\minipr_h\}$ and $\nonm=\pri_{\omega,h}\leq\max\{\dd,\supipr_h\}$. 
		\end{itemize}
	\end{lem}

	\section{connections between localizations and slalom predictions} \label{sec_connection}
	
	In this section we see the connection between localizations and slalom predictions.
	First, predictions are easier than localizations because of the hints of initial segments of functions, so we immediately have:
	\begin{lem}
		\label{pr_lq_lc}
		$\gprbh\lq\glcbh$, $\lprbh\lq\llcbh$ and $\iprbh\lq\ilcbh$ for any $b$ and $h$.
	\end{lem}
	
	However, regardless of such hints, predictions can be harder than localizations when the width of slaloms of localizations is sufficiently wider:

	\begin{dfn}
		For $h\in\oo$, define $h^\prime_{h}\in\oo$ by $h^\prime_{h}(n)\coloneqq n\cdot\prod_{i\leq n}h(i)$.
	\end{dfn}
	
	\begin{lem}(\cite{Bla10})
		\label{lem_Bla10}
		Let $b\in(\omega+1)^\omega$ and $h\in\oo$ such that $b>h^\prime\coloneqq h^\prime_{h}$, where $<$ denotes the total strict domination on $(\omega+1)^\omega$. Then, $\mathbf{GLc}(b,h^\prime)\lq\gprbh$, $\mathbf{LLc}(b,h^\prime)\lq\lprbh$ and $\mathbf{ILc}(b,h^\prime)\lq\iprbh$.
	\end{lem}
	
	\begin{proof}
		Let $\pi$ be a global $(b,h)$-predictor. 
		For $s\in\seqb$, take a tree $T^s\subseteq\seqb$ of stem $s$ which contains all $t\in\seqb$ such that $s\subseteq t$ and $t(n)\in\pi(t\on n)$ for all $n\in|t|\setminus |s|$ and we may assume that every $t\in T^s\cap\omega^{n}$ is $h(n)$-branching for all $n\geq |s|$. Put $\pi_s(n)\coloneqq\bigcup\{t(n):t\in T^s\cap \omega^{n+1}\}$ for $n<\omega$ and note that $|\pi_s(n)|\leq \prod_{i\leq n}h(i)$. Enumerate $\{s_i:i<\omega\}=\seqb$ and let $\varphi_\pi(n)\coloneqq\bigcup\{\pi_{s_i}(n):i< n\}$ for $n<\omega$ and note that $|\varphi_\pi(n)|\leq \prod_{i\leq n}h(i)\cdot n=h^\prime_{h}(n)$ , so $\varphi_\pi$ is a $(b,h^\prime_{h})$-slalom.
		Take $x\in \prod b$ and $n_0<\omega$ arbitrarily and put $s_{i}\coloneq x\on n_0$.
		Then, for any $n>i$, $x(n)\in\pi(x\on n)$ implies $x(n)\in\varphi_\pi(n)$, which induces the three Tukey connections. 
	\end{proof}

	Consequently, localizations and predictions have the same limit value in the following sense: 
	\begin{thm}
		\label{thm_limit_same_pre}
		Let $h\in\oo$.
		\begin{enumerate}
			\item $\mingpr_h\leq \minglc$ and $\supgpr_h\geq\supglc$.
			In particular, $\mingpr_{h}= \minglc$ and $\supgpr_{h}=\supglc$ if $h$ goes to infinity by Lemma \ref{pr_lq_lc} and Theorem \ref{thm_GILc_h_indep}. 
			\item(e.g. \cite[Proposition 2.2]{Laf97}) $(\ee_{ubd}=)\minlpr_h=\minllc$ and  $(\pr_{ubd}=)\suplpr_h=\supllc$.
			\item $\minipr_h= \minilc$ and $\supipr_h=\supilc$.
		\end{enumerate}
	\end{thm}
	
	\begin{proof}
		We only prove the former (in)equalities.
		\begin{enumerate}
			\item 
			Since $h^\prime\coloneqq h^\prime_h$ goes to infinity, $\minglc=\minglc_{h^\prime}=\min\{\bbglc_{b,h^\prime}: b\in\oo\}=\min\{\bbglc_{b,h^\prime}: b\in\oo, b>h^\prime\}$ holds. Let $b\in\oo$ be arbitrary with $b>h^\prime$. By Lemma \ref{lem_Bla10}, $\mingpr_h\leq \eeg_{b,h}\leq \bbglc_{b,h^\prime}$. Since $b>h^\prime$ was arbitrary, we have $\mingpr_h\leq \minglc$.
			\item If $h$ goes to infinity, $\minlpr_h=\minllc$ is obtained similarly, so we may assume $h$ is bounded. 
			Since $\minlpr_h\leq \minllc$ is proved similarly as well, we show the converse. Let $b\in\oo$ be arbitrary, $h^+\geq h$ be any function going to infinity and $b^+$ witness $\mathbf{LPr}(b,h)\lq\mathbf{LPr}(b^+,h^+)$ by Lemma \ref{lem_their_3.14_pred_1}.
			Since $h^+$ goes to infinity, $\minllc=\minlpr_{h^+}\leq\eel_{b^+,h^+}\leq\eel_{b,h}$.
			Since $b$ was arbitrary, we obtain $\minllc\leq\minlpr_{h}$. 
			\item $\minipr_h\leq \minilc$ is proved similarly and $\minipr_h\geq\minilc_h=\minilc$ by Lemma \ref{pr_lq_lc} and Theorem \ref{thm_GILc_h_indep}.
		\end{enumerate}
	\end{proof}
	
	\begin{rem}
		If $h$ is bounded, $\mingpr_h< \minglc$ might happen: 
		It is not hard to see that $\sigma$-$(k+1)$-linked forcings are $\mathbf{GPr}(k+1,k)$-good and hence $\mingpr_k<\addn=\bbglc_{\omega,h}\leq\minglc$ holds in the Amoeba model (see also Theorem \ref{thm_bbglcoh_addn}, \cite[Lemma 3.1, Theorem 3.8(a)]{BS01}).
	\end{rem}

	In the local and infinite versions, the same phenomenon that predictions can be harder than localizations happens in another case, namely, when the whole space of predictions is sufficiently larger:
	\begin{thm}
		\label{thm_connection}
		Let $b, h, h^-\in\oo$ and assume $h$ goes to infinity.
		Then, there is $b^+\in\oo$ such that $\mathbf{LLc}(b,h)\lq\mathbf{LPr}(b^+,h^-)$ and $\mathbf{ILc}(b,h)\lq\mathbf{IPr}(b^+,h^-)$.
	\end{thm}
	
	\begin{proof}
		Inductively take natural numbers $0=i_{-1}<i_0<\cdots$ which satisfy for all $k<\omega$,
		\begin{equation*}
			\prod_{i<i_{k-1}}b(i)\cdot h^-(k)\leq h(i_k-1),
		\end{equation*}
		which is possible since $h$ goes to infinity.
		Define $b^+$ by $b^+(k)\coloneq\prod_{i\in I_k}b(i)$ where $I_k\coloneq\left[i_{k-1},i_k \right)$ for each $k<\omega$.
		Assume we are given a global $(b^+,h^-)$-predictor $\pi$. For $k<\omega$, define:
		\begin{equation*}
			A_k\coloneq\left\{t(i_k-1):t\in\pi(\langle\sigma\on I_{k^\prime}:k^\prime<k\rangle),~\sigma\in\prod_{i<i_{k-1}}b(i)\right\},
		\end{equation*}
		which is valid since if $\sigma\in\prod_{i<i_{k-1}}b(i)$, by identifying functions and natural numbers, $\langle\sigma\on I_{k^\prime}:k^\prime<k\rangle\in\prod_{k^\prime<k}b^+(k^\prime)$ and if $t\in\pi(\langle\sigma\on I_{k^\prime}:k^\prime<k\rangle)$, $t$ belongs to $b^+(k)=\prod_{i\in I_k}b(i)$ and again by identification $t$ is a function and $t(i_k-1)\in b(i_k-1)$.
		Also note that $|A_k|\leq\prod_{i<i_{k-1}}b(i)\cdot h^-(k)\leq h(i_k-1)$.
		Therefore, $A_k\in[b(i_k-1)]^{\leq h(i_k-1)}$, so there is a $(b,h)$-slalom $\varphi_\pi$ satisfying $\varphi_\pi(i_k-1)=A_k$ for all $k<\omega$.
		By construction, for any $k<\omega$ and $f\in\prod b$,
		\begin{equation*}
			f\on I_k\in\pi(\langle f\on I_{k^\prime}:k^\prime<k\rangle)\text{ implies } f(i_k-1)\in\varphi_\pi(i_k-1),
		\end{equation*}
		which induces $\mathbf{LLc}(b,h)\lq\mathbf{LPr}(b^+,h^-)$ and $\mathbf{ILc}(b,h)\lq\mathbf{IPr}(b^+,h^-)$.
	\end{proof}
	
	\section{Connection with ideals} \label{sec_ideals}
	There are many known results on the connections between ideals on the reals and global/infinite localization cardinals. Here are such examples:
	\begin{thm}(\cite{Bar84}, \cite[Theorem 4.2]{CM23})
		\label{thm_bbglcoh_addn}
		$\bbglc_{\omega,h}=\addn$ and $\ddglc_{\omega,h}=\cofn$ for $h\in\oo$ going to infinity.
	\end{thm}
	
	
	\begin{thm}(\cite{Mil82}, \cite{Bar87}, \cite[Theorem 5.1]{CM23})
		\label{thm_ddilc_nonm}
		$\ddilc_{\omega,h}=\nonm$ and $\bbilc_{\omega,h}=\covm$ for $h\in\oo$.
	\end{thm}
	
	\begin{thm}(\cite[Lemma 2.3]{KM22},\cite[Lemma 6.2]{CM23})
		\label{thm_N_E}
		\begin{enumerate}
			\item If $\sum_{n<\omega}\frac{h(n)}{b(n)}<\infty$, then $\covn\leq \ddilc_{b,h}$ and $\bbilc_{b,h}\leq\nonn$. \label{item_N_E_N}
			\item If $\sum_{n<\omega}\frac{h(n)}{b(n)}=\infty$, then $\cov(\mathcal{E})\leq \bbilc_{b,h}$ and $\ddilc_{b,h}\leq\none$, where $\cal{E}$ denotes the $\sigma$-ideal generated by closed null sets.
		\end{enumerate}
	\end{thm}
	
	Now let us look at local cases. First, we easily have the following by considering the set of functions predicted by a single local predictor: 
	\begin{lem}
		\begin{itemize}
			\item Let $b\in(\omega+1)^\omega$ and $h\in\oo$. Then, 
			$\bbllc_{b,h}\leq\eel_{b,h}\leq\nonm$ and $\cov(\cal{M})\leq\prl_{b,h}\leq\ddllc_{b,h}$.
			\item Let $b,h\in\oo$. If $\limsup_n\frac{h(n)}{b(n)}<1$, then 
			$\bbllc_{b,h}\leq\eel_{b,h}\leq\none$ and $\cov(\cal{E})\leq\prl_{b,h}\leq\ddllc_{b,h}$.  
		\end{itemize}
	\end{lem}
	
	\begin{rem}
		If $\limsup_n\frac{h(n)}{b(n)}=1$, $\none<\bbllc_{b,h}$ and $\ddllc_{b,h}<\cov(\cal{E})$ might happen: The author introduced in \cite{Yam24calE} the forcing notion $\lebb_b$ which adds an $\mathbf{LLc}(b,b-1)$-unbounded real and proved that $\lebb_b$ keeps $\none$ small if $b(n)\geq 2^n$ for $n<\omega$. Consequently, in his model constructed in \cite[Theorem 5.9]{Yam24calE}, $\none<\bbllc_{b,b-1}=\nonm<\covm=\ddllc_{b,b-1}<\cov(\cal{E})$ holds.
	\end{rem}

	Let us focus on the uniformity of the null and meager additive ideals $\mathcal{NA}$ and $\mathcal{MA}$.
	
	\begin{dfn}
		Let $\mathcal{I}$ be an ideal on $2^\omega$. $\mathcal{IA}$ denotes the set of all $\mathcal{I}$-additive sets, i.e., the collection of $X\subseteq2^\omega$ such that $A+X\in\mathcal{I}$ for all $A\in\mathcal{I}$. Here, the addition $+$ on $2^\omega$ is defined by identifying $2^\omega\cong (\mathbb{Z}/2\mathbb{Z})^\omega$.  
	\end{dfn}
	
	The additivity and uniformity of $\mathcal{NA}$ can be characterized using the limit global number:
	
	\begin{fac}(\cite[Lemma 2.2]{Paw85}, \cite[Corollary 1.7, Theorem A]{CMR24})
		$\add(\mathcal{NA})=\non(\mathcal{NA})=\minglc$.
	\end{fac}
	
	In particular, $\non(\mathcal{NA})=\minglc\leq\minllc=\minlpr=\ee_{ubd}$ holds, which is already shown by Cardona in \cite[Theorem 2.1]{Car24}. However, in the case of $\mathcal{MA}$, we show that the opposite direction holds:
	
	\begin{thm}
		\label{thm_eubd_nonMA}
		$\ee_{ubd}\leq\non(\mathcal{MA})$.
	\end{thm}
	
	To prove this theorem, let us introduce the following relational system to characterize the uniformity of $\mathcal{MA}$:
	
	\begin{dfn}(\cite[Definition 2.7]{CMR24})
		Let $\ip$ denote the set of all interval partitions of $\omega$.
		\begin{enumerate}[label=(\arabic*)]
			\item  For $f, g\in\oo$ and $I=\langle I_k:k<\omega\rangle\in\ip$, $f  \sqsubset^\bullet (I,g)$ if for all but finitely many $k<\omega$, there exists $i\in I_k$ such that $f(k)=g(k)$.
			\item  For $b\in\oo$, denote the relational system $\R_b\coloneq\langle\prod b,\ip\times\prod b,\sqsubset^\bullet\rangle$.
		\end{enumerate}

	\end{dfn} 
	
	\begin{fac}(\cite[Theorem 2.2]{BJ94}, \cite[Lemma 2.10]{CMR24})
		\label{fac_chara_nonMA}
		$\non(\mathcal{MA})=\min\{\bb(\R_b):b\in \oo\}$.
	\end{fac}
	
	$\R_{b}$ is Tukey below some $\mathbf{LLc}(b^+,h)$:
	
	\begin{lem}
		\label{lem_con_llc_and_Rb}
		Let $b,h\in\oo$. Then, there exists $b^+\in\oo$ such that $\R_b\lq\mathbf{LLc}(b^+,h)$.
	\end{lem}
	
	\begin{proof}
		Take some $I=\langle I_k:k<\omega\rangle\in\ip$ such that $|I_k|\geq h(k)$.
		Define $b^+(k)\coloneq \prod_{i\in I_k}b(i)$.
		Let $\varphi$ be any $(b^+,h)$-slalom and fix $k<\omega$.
		Let $S\coloneq\varphi(k)\in b^+(k)\coloneq \prod_{i\in I_k}b(i)$, which is assumed to be a set of functions on $I_k$.
		Enumerate $\{ t_i:i<|S|\}=S$ and let $I^-\coloneq\{\min I_k+i:i<|S|\}\subseteq I_k$, since $|S|\leq h(k)\leq|I_k|$.
		define a function $g_k$ on $I^-$ by:
		\begin{equation*}
			g_k(\min I_k+i)\coloneq t_i(i)\in b(i).
		\end{equation*}
		Unfix $k$ and by collecting all $g_k$ together, let $g=g_\varphi\in\prod b$ be some function which satisfies $g(\min I_k+i)=g_k(\min I_k+i)$ for all $k<\omega$ and $i<|\varphi(k)|$.
		By construction, note that for any $x\in \prod b$ and $k<\omega$,
		\begin{equation*}
			x \on I_k\in\varphi(k)\text{ implies }x(i)=g_\varphi(i)\text{ for some }i\in I_k.
		\end{equation*}
		Let $D\in\ooo$ and $I_D=\langle I^\prime_k:k<\omega\rangle\in\ip$ be such that if $k$ is the $j$-th element of $D$, then $I_k\subseteq I^\prime_j$.
		Again by construction, for $x\in\prod b$, $\langle x\on I_k:k<\omega\rangle\in^*(D,\varphi)$ implies $x\sqsubset^\bullet(I_D,g_\varphi)$, which shows $\R_b\lq\mathbf{LLc}(b^+,h)$.
	\end{proof}
	
	\begin{proof}[Proof of Theorem \ref{thm_eubd_nonMA}]
		Let $h$ be any function going to infinity.
		By Theorem \ref{thm_limit_same_pre},
		we have $\ee_{ubd}=\minlpr_h=\minllc_h$, so we show $\minllc_h\leq\non(\mathcal{MA})$ instead.
		Let $b\in \oo$ be arbitrary. By Lemma \ref{lem_con_llc_and_Rb}, some $b^+\in\oo$ satisfies $\R_b\lq\mathbf{LLc}(b^+,h)$, so we have $\minllc_h\leq\bbllc_{b^+,h}\leq\bb(\R_b)$. Since $b$ was arbitrary, by Fact \ref{fac_chara_nonMA} we have $\minllc_h\leq \min\{\bb(\R_b):b\in \oo\}=\non(\mathcal{MA})$.
	\end{proof}
	
	\begin{rem}
		$\ee_{ubd}<\non(\mathcal{MA})$ holds in the Hechler model: Brendle (essentially) showed in \cite{Bre95} that $\ee_{ubd}\leq\ee_{\mathrm{id}}=\aleph_1$ in the Hechler model ($\mathrm{id}$ denotes the identity function on $\omega$), and in that model $\non(\mathcal{MA})\geq\addm=\min\{\bb,\covm\}=2^{\aleph_0}>\aleph_1$ holds. 
	\end{rem}
	
	Figure \ref{fig_numbers} illustrates the relationship of the cardinal invariants below $\nonm$ we have seen and here are additional explanations:
	\begin{itemize}
		\item $\covn\leq\mathrm{supI}$ by Theorem \ref{thm_N_E}\eqref{item_N_E_N}.
		\item $\non(\mathcal{MA})\leq\none$ by e.g. \cite[Corollary 1.12]{CMR24}.
		\item $\add(\mathcal{M})\leq\non(\mathcal{MA})$ by $\add(\mathcal{M})\leq\add(\mathcal{MA})\leq\non(\mathcal{MA})$.
		\item $\min\{\ee,\bb\}\leq\addm$ by $\addm=\min\{\bb,\covm\}$ and $\ee\leq\covm$.
		\item $\ddilc_{b,h}$ can be finite by \cite[Theorem 3.13]{CM23}.
	\end{itemize}
	
	\begin{figure}[h]
		\centering 
		\begin{tikzpicture}
			\tikzset{
				textnode/.style={text=black}, 
			}
			\tikzset{
				edge/.style={color=black, thin}, 
			}
			\tikzset{cross/.style={preaction={-,draw=white,line width=9pt}}}
			\newcommand{\w}{2.8}
			\newcommand{\h}{3.3}
			\newcommand{\eqsc}{0.3}
			\newcommand{\mx}{4.6*\w}
			\newcommand{\lx}{3.7*\w}
			\newcommand{\rx}{5.3*\w}
			\newcommand{\dx}{0.4*\w}
			
			\newcommand{\llx}{3.25*\w}
			
			\node[textnode] (addN) at (3*\w,  0) {$\addn$};
			\node[textnode] (covN) at (3*\w,  \h*2) {$\covn$};

			\node[textnode] (addM) at (\w*6,  0) {$\addm$};
			\node[textnode] (b) at (\w*6,  \h) {$\bb$};
			\node[textnode] (nonM) at (\w*6,  \h*2) {$\nonm$};
			
			\node[textnode] (covM) at (\w*6.7,  0) {$\covm$};


			\node[textnode] (bbglcoh) at (3*\w,  -\eqsc*\h) {$\bbglc_{\omega,h}$};
			\node[textnode] (eegoh) at (3*\w,  -2*\eqsc*\h) {$\eeg_{\omega,h}$};

			\node[textnode] (mineb) at (\mx,  0) {$\min\{\ee,\bb\}$};
			\node[textnode] (bllcoh) at (\mx,   -\eqsc*\h) {$\bbllc_{\omega,h}$};
			
			\node[textnode] (e) at (\mx,  0.3*\h) {$\ee$};
			\node[textnode] (eloh) at (\mx,  0.6*\h) {$\eel_{\omega,h}$};
			\node[textnode] (eubd) at (\mx,  0.9*\h) {$\mathrm{minL}=\ee_{ubd}$};
			
			\node[textnode] (minglc) at (\lx,  0.3*\h) {$\mathrm{minG}=\non(\mathcal{NA})$};

			
			\node[textnode] (supI) at (\mx,  2*\h) {$\mathrm{supI}$};
			\node[textnode] (nonMA) at (\rx-0.05*\w,  1.25*\h) {$\non(\mathcal{MA})$};
			\node[textnode] (nonE) at (\rx+0.3*\w,  1.6*\h) {$\none$}; 
			
			\node[textnode] (bllcbh) at (\mx,  1.25*\h) {$\bbllc_{b,h}$}; 
			\node[textnode] (elbh) at (\mx,  1.6*\h) {$\eel_{b,h}$}; 
			
			\node[textnode] (bglcbh) at (\lx,  1.25*\h) {$\bbglc_{b,h}$}; 
			\node[textnode] (egbh) at (\lx,  1.6*\h) {$\eeg_{b,h}$}; 
			
			\node[textnode] (dilcbh) at (\llx,  1.7*\h) {$\ddilc_{b,h}$}; 
			\node[textnode] (pribh) at (\llx,  1.3*\h) {$\pri_{b,h}$};

			\node[textnode] (aleph1) at (2.5*\w,  0) {$\aleph_1$};
			\draw[->, edge] (aleph1) to (addN);

			\draw[->, edge] (addN) to (covN);
			\draw[->, edge] (addN) to (mineb);
			\draw[->, edge] (mineb) to (addM);
			
			\draw[->, edge] (covN) to (supI);
			\draw[->, edge] (supI) to (nonM);

			\draw[->, edge] (addM) to (b);
			\draw[->, edge] (b) to (nonM);
			\draw[->, edge] (addM) to (covM);
			
			
			\draw[->, edge] (addN) to (minglc);
			\draw[->, edge] (minglc) to (eubd);	
			\draw[->, edge] (mineb) to (e);
			\draw[->, edge] (e) to (eloh);
			\draw[->, edge] (eloh) to (covM);
			\draw[->, edge] (eloh) to (eubd);
			
			\draw[->, edge] (eubd) to (nonMA);
			\draw[->, edge] (nonMA) to (nonE);
			\draw[->, edge] (nonE) to (nonM);
			

			\draw[->, edge] (eubd) to (bllcbh);
			\draw[->, edge] (bllcbh) to (elbh);
			\draw[->, edge] (elbh) to (nonM);
			
			\draw[->, edge] (bglcbh) to (egbh);
			\draw[->, edge] (minglc) to (bglcbh);
			\draw[->, edge] (bglcbh) to (egbh);
			\draw[->, edge] (egbh) to (elbh);
			\draw[->, edge] (bglcbh) to (bllcbh);
			
			\draw[->, edge] (pribh) to (dilcbh);
			\draw[->, edge] (dilcbh) to (supI);

			\draw[double distance=2pt, edge] (addN) to (bbglcoh);
			\draw[double distance=2pt, edge] (bbglcoh) to (eegoh);
			\draw[double distance=2pt, edge] (mineb) to (bllcoh);
			
			\draw[->, edge] (addM) to (nonMA);

		\end{tikzpicture}
		\caption{Diagram of cardinal invariants below $\nonm$. Here, $b\in(\omega\setminus2)^\omega$ and $h\in\prod b$ go to infinity and $\mathrm{minG}\coloneqq\minglc=\mingpr_h$, $\mathrm{minL}\coloneqq\minllc=\minlpr_{h^\prime}$, $\mathrm{supI}\coloneqq\supilc=\supipr_{h^\prime}$ for $h^\prime\geq 1$. }\label{fig_numbers}
	\end{figure}

	\begin{acknowledgements}
		This note was written for the proceedings of the RIMS Set Theory Workshop 2024 \textit{Recent Developments in Axiomatic Set Theory}, held at Kyoto University RIMS. The author is grateful to the organizer Masahiro Shioya for letting him give a talk at the Workshop and submit an article to the proceedings.
		The author also thanks his supervisor J\"{o}rg Brendle for his helpful comments. 
		This work was supported by JST SPRING, Japan Grant Number JPMJSP2148.
	\end{acknowledgements}
	
	\bibliographystyle{alpha}
	\bibliography{../ref.bib}

\end{document}